\documentclass [12pt,a4paper]{article}
\usepackage[T1]{fontenc}
\usepackage[ansinew]{inputenc}
\usepackage{amssymb}
\usepackage{amsmath}
\usepackage{graphicx}
\usepackage{amsthm}
\usepackage{cite}

\newtheorem{theorem}{\textbf{Theorem}}[section]
\newtheorem{lemma}{\textbf{Lemma}}[section]

\newtheorem{remark}{\textbf{Remark}}[section]
\newtheorem{definition}{\textbf{Definition}}[section]

\allowdisplaybreaks[4]

\def\be{\begin{equation}}
\def\ee{\end{equation}}
\def\bea{\begin{eqnarray}}
\def\eea{\end{eqnarray}}
\def\bt{\begin{theorem}}
\def\et{\end{theorem}}
\def\bl{\begin{lemma}}
\def\el{\end{lemma}}
\def\br{\begin{remark}}
\def\er{\end{remark}}
\def\bc{\begin{corollary}}
\def\ec{\end{corollary}}
\def\bd{\begin{definition}}
\def\ed{\end{definition}}

\def\s0t{\sup _{0\leq \tau\leq t}}
\def\C0T{C([0,T];\,}

\def\DAS{D( A^{\frac{S}{2}})}
\def\DAS1{D( A^{ \frac{S+1}{2}})}

\def\p{\partial}

\def\dis{\displaystyle}

 \def\non{\nonumber }

\begin{document}
\noindent \underline{Title}: Convergence to Equilibrium for the Cahn-Hilliard Equation with Wentzell Boundary Condition\\[5mm]
\underline{Author}:\\
\quad Hao {\sc Wu}\\
School of Mathematical Sciences\\
Fudan University\\
Shanghai 200433, P.R. China\\[3mm]

\noindent \underline{Mailing address}:\\
Hao Wu\\
School of Mathematical Sciences\\
Fudan University\\
Han Dan Road No. 220, Shanghai, 200433,\\
P.R. China\\
Email: haowufd@yahoo.com \\
\bigskip


\newpage
\setcounter{page}{1}

\title{Convergence to Equilibrium for the Cahn-Hilliard Equation with Wentzell Boundary Condition}

\author{ Hao Wu  \\ School of Mathematical Sciences\\ Fudan University \\
 200433 Shanghai, P.R. China\\  haowufd@yahoo.com }

\date{\today}

\maketitle


\begin{abstract}
In this paper we consider the Cahn-Hilliard equation endowed with
Wentzell boundary condition which is a model of phase separation in
a binary mixture contained in a bounded domain with permeable wall.
Under the assumption that the nonlinearity is analytic with respect
to unknown dependent function, we prove the convergence of a global
solution to an equilibrium as time goes to infinity by means
of a suitable \L ojasiewicz-Simon type inequality with boundary term. Estimates of convergence rate are also provided.\\
\textbf{Keywords}: Cahn-Hilliard equation, Wentzell boundary
conditions, \L ojasiewicz-Simon inequality, convergence to
equilibrium.

\end{abstract}


\section{Introduction}
\setcounter{equation}{0}

This paper is concerned with the asymptotic behavior of the global solution to the following Cahn-Hilliard equation
\begin{equation}
u_t=\Delta\mu,\qquad \text{in} \ \ [0,T]\times\Omega\label{cw1}
\end{equation}
\begin{equation}
\mu=-\Delta u+f(u),\qquad \text{in} \ \
[0,T]\times\Omega,\label{cw2}
\end{equation}
subject to Wentzell boundary condition
\begin{equation}
\Delta\mu+b \partial _\nu \mu+c\mu=0,\qquad \text{on} \ \
[0,T]\times\Gamma,\label{cw3}
\end{equation}
the variational boundary condition
\begin{equation}
-\alpha\Delta_\parallel u+\partial _\nu u+\beta
u=\frac{\mu}{b},\qquad \text{on} \ \ [0,T]\times\Gamma, \label{cw4}
\end{equation}
and initial datum
\begin{equation}
u(0,x)=\psi_0, \qquad \text{in}\ \  \Omega.\label{cw5}
\end{equation}
In above, $0<T\leq \infty$, $\Omega $ is a bounded domain in
$\mathbb{R}^n\ (n=2,3)$ with smooth boundary $\Gamma$.
$\alpha,\beta,b,c$ are positive constants. $\Delta _\parallel $ is
the Laplace-Beltrami operator on $\Gamma$,
 and $\nu$ is the outward normal direction to the boundary.

The Cahn-Hilliard equation arises from the study of spinodal
decomposition of binary mixtures that appears, for example, in
cooling process of alloys, glass or polymer mixtures (see
\cite{CH58,K01,NS84,T} and the references cited therein). $\mu$ is
called chemical potential in the literature. The classical
Cahn-Hilliard equation is equipped with the following homogeneous
Neumann boundary conditions
 \bea \p
_\nu \mu &=&0, \qquad t>0,\ x
\in \Gamma,\label{1.2a}\\
 \partial_\nu u&=&0, \qquad t > 0, \ x \in \Gamma.\label{1.2aa}
  \eea
Boundary \eqref{1.2a} has a clear physical meaning: there cannot be
any exchange of the mixture constituents through the boundary
$\Gamma$ which implies that the total mass $\int_\Omega u dx$ is
conversed for all time. The boundary condition \eqref{1.2aa} is
usually called variational boundary condition which together with
\eqref{1.2a} result in decreasing of the following bulk free energy
\be E_b(u)=\int_\Omega\left(\frac12|\nabla u|^2+F(u)\right)dx,\ee
where $F(s)=\int^s_0f(z)dz$. A typical example in physics for
potential $F$ is the so-called 'double-well' potential
$F(u)=\frac14(u^2-1)^2$.

 For the equations \eqref{cw1} \eqref{cw2} subject to boundary conditions \eqref{1.2a} \eqref{1.2aa}
  and initial datum \eqref{cw5}, extensive study
has been made. We refer e.g., to \cite{EZ86,Z86,NP93,GN95,T,RH} and
the references cited therein. In particular, convergence to
equilibrium for the global solution in higher space dimension case
was proved in \cite{RH}.

Recently, a new model has been derived when the effective
interaction between the wall (i.e., the boundary $\Gamma$) and two
mixture components are short-ranged (see Kenzler et al. \cite{K01}).
In such a situation, it is pointed out in \cite{K01} that, the
following surface energy functional
 \be
E_s(u)=\int_\Gamma \left(\frac{\sigma_s}{2}\left|\nabla_\parallel u
\right|^2 +
 \frac{g_s}{2} u^2 - h_s u\right) dS,\label{cese}\ee
 with $\nabla_\parallel$ being the covariant gradient operator on
 $\Gamma$ (see e.g. \cite{MZ05}), should be added to the bulk free energy $E_b(u)$ to form a total
 free
energy functional \be E(u)=E_b(u)+E_s(u).\ee In above, $\sigma_s>0$,
$g_s>0$, $h_s\neq 0$ are given constants. Together with the no-flux
boundary \eqref{1.2a} condition, the following dynamical boundary
condition is posed in order that the total energy $E(u)$ is
decreasing with respect to time:
 \be \sigma_s \Delta _\parallel u -
\partial_\nu u - g_s u+ h_s =\displaystyle \frac{1}{\Gamma_s} u_t,\qquad t>0,\ x
\in \Gamma.\label{dyna}
 \ee
We refer to \cite{CEP06,RZ03, PRZ,WZ1,MZ05} for extensive study for
system \eqref{cw1}\eqref{cw2} with boundary conditions
\eqref{1.2a}\eqref{dyna} and initial datum \eqref{cw5}. In
particular,  Wu $\&$ Zheng \cite{WZ1} proved the convergence to
equilibrium for a global solution as time goes to infinity by
deriving a new type of \L ojasiewicz-Simon inequality with boundary
term (see also \cite{CEP06} for a different proof).

Based on the above model, in a quite recent article by Gal
\cite{Gal}, the author proposed \eqref{cw1}--\eqref{cw5} as a
variation model which describes phase separation in a binary mixture
confined to a bounded region $\Omega $ with porous walls. Instead of
the no-flux boundary condition \eqref{1.2a}, the Wentzell boundary
condition \eqref{cw3} is derived from mass conservation laws that
include an external mass source (energy density) on boundary
$\Gamma$. This may be realized, for example, by an appropriate
choice of the surface material of the wall, i.e., the wall $\Gamma $
may be replaced by a penetrable permeable membrane (ref.
\cite{Gal}). Then, the variational boundary condition \eqref{cw4} is
introduced in order that the system \eqref{cw1}--\eqref{cw5} tends
to minimize its total energy
 \be
 E(u)= \int_{\Omega} \left(\frac{1}{2} |\nabla u|^2 + F(u)\right) dx
 + \int_\Gamma\left(\frac{\alpha}{2}\left|\nabla_\parallel u \right|^2 +\frac{\beta}{2} u^2\right) dS. \ee
 Namely,
 \be\frac{d}{dt}E(u(t))=-\int_\Omega |\nabla \mu|^2dx
 -\frac{c}{b}\int_\Gamma \mu^2 dS\leq 0.\ee
 For
more intensive discussions, we refer to \cite{Gal, Gal2, Gold06}.

In \cite{Gal}, the existence and uniqueness of global solution to
problem \eqref{cw1}-\eqref{cw5} has been proved by adapting the
approach in \cite{RZ03}. Later, in \cite{Gal2}, the same author
studied the problem in a further way that he obtained the existence
and uniqueness of a global solution to the problem under more
general assumptions than those in \cite{Gal}. He showed that the
global solution defines a semiflow on certain function spaces and
also proved the existence of an exponential attractor with finite
dimension.

Then a natural question is: whether the global solution of system
\eqref{cw1}-\eqref{cw5} will converge to an equilibrium as time goes
to infinity? This is just the main goal of this paper. Moreover, we
shall provide estimates for the rate of the convergence (in  higher
order norm).

\br\label{r1} Without loss of generality, in the following text, we
set positive constants $b,c, \alpha, \beta$ to be $1$.  In this
paper, we simply use $\|\cdot\|$ for the norm on $L^2(\Omega)$ and
equip $H^1(\Omega)$ with the equivalent norm
 \be
  \|u\|_{H^1(\Omega)}=\left(\int_\Omega|\nabla u|^2dx+ \int_\Gamma
 u^2dS\right)^{1/2}.\label{cwH1}\ee\er

Before stating our main result, first we make some assumptions
on nonlinearity $f$. \\
(\textbf{F1}) $f(s)$ is analytic in $s\in \mathbb{R}$.\\
(\textbf{F2})
$$|f(s)|\leq C(1+|s|^p),\;\;\;\; \forall
  s\in \mathbb{R},$$
 where $C\geq 0$, $p>0$ and $p\in (0,5)$ for $n=3$.\\
(\textbf{F3})
 \be \liminf_{|s|\rightarrow \infty} f'(s)>0.\non\ee
 \br Assumption (\textbf{F1}) is made so that we are able to derive an
extended \L ojasiewicz-Simon inequality to prove our convergence
result. Assumption (\textbf{F2})  implies that the nonlinear term
has a subcritical growth. Assumption (\textbf{F3}) is some kind of
dissipative condition. (\textbf{F3}) is supposed in \cite{Gal, Gal2}
to obtain the existence and uniqueness of global solution to the
evolution problem \eqref{cw1}--\eqref{cw5}. Moreover, (\textbf{F3})
together with (\textbf{F2}) enable us to prove the existence result
for stationary problem \eqref{cwsta} by variational method (see
Section 3). It's easy to check that the nonlinearity $f(u)=u^3-u$
corresponding to the most important physical potential
$F(u)=\frac14(u^2-1)^2$ satisfies all the assumptions stated above.
\er
 Let $V$ be the Hilbert space which, as
introduced in \cite{Gal}, is the completion of $C^1(\Omega)$ with
the following inner product and the associated norm:
 \be (u, v)_V=\int_{\Omega} \nabla u \cdot
\nabla v dx +\int_{\Gamma}\left(\nabla_{\parallel}u \cdot
\nabla_{\parallel}v + u v\right) d S,\quad \forall\ u, v \in
V.\label{cw1.4}\ee

The main result of this paper is as follows.
 \bt \label{cwcon}  Let (\textbf{F1})--(\textbf{F3}) be satisfied.
 For any initial
datum $u_0 \in V$, the solution $u(t,x)$ to problem
(\ref{cw1})--(\ref{cw5}) converges to a certain equilibrium
$\psi(x)$ in the topology of $H^3(\Omega)\cap H^3(\Gamma)$ as time
goes to infinity, i.e., \be \dis{\lim_{t\rightarrow
+\infty}}\left(\|u(t,\cdot)-\psi\|_{H^3(\Omega)}+\|u(t,
\cdot)-\psi\|_{H^3(\Gamma)}\right)=0. \label{cwconr}\ee Moreover, we
have the following estimate for the rate of convergence: \be
\|u-\psi\|_{H^3(\Omega)}+\|u-\psi\|_{H^3(\Gamma)}+\|u_t\|_{V}\leq
 C(1+t)^{-\theta/(1-2\theta)},\qquad t\geq 1.\label{cwconrate}\ee
Here, $C\geq 0$, $\psi(x)$ is an equilibrium to problem
(\ref{cw1})--(\ref{cw5}), i.e., a solution to the following
nonlinear boundary value
problem: \be \left\{\begin{array}{l}  - \Delta \psi +f(\psi)=0, \;\;\;\; x \in \Omega,\\
-\Delta_\parallel \psi+\partial _\nu \psi +\psi =0, \;\;\; x\in
\Gamma,
  \end{array}
 \label{cwsta}
 \right.
 \ee
and $\theta\in (0,\frac12)$ is a constant depending on $\psi(x)$.
\et

\par Before giving the detailed proof of Theorem 1.1, let's first
  recall some related results
  in the literature.
  The study of asymptotic behavior of solutions
  to nonlinear dissipative evolution equations has attracted a
  lot of interests of many mathematicians for a long period of time.
  Unlike in 1-d case (see \cite{M78,Ze}), the situation in higher space dimension case can be quite
  complicated. On one hand, the topology of the set of stationary solutions can be non-trivial and may
   form a continuum. On the other hand, a counterexample has been given in \cite{PS}
  for a semilinear parabolic equation saying that even the nonlinear term being $C^\infty$
  cannot ensure the convergence to a single equilibrium (see also \cite{PR96}).
  In 1983 Simon in \cite{S83} proved that for a semilinear parabolic equation if the
  nonlinearity is analytic in unknown function $u$, then
  convergence to equilibrium for bounded global solutions holds. His idea relies on generalization of the
  \L ojasiewicz inequality (see \cite{L1,L2,L3}) for analytic functions defined in finite
  dimensional space $\mathbb{R}^m$. Since then, Simon's idea has been applied to prove convergence results for
  many evolution equations, see e.g., \cite{J981,J982,HT01,HJ01,RH} and the references cited
  therein. To the best of our knowledge, most previous work are
  concerned with evolution equations subject to homogeneous Dirichlet or Neumann boundary
  conditions.

   Our problem (\ref{cw1})--(\ref{cw5}) has the following features. The first boundary condition (\ref{cw3}) is Wentzell boundary condition which involves the time derivative of $u$;
   the second boundary condition (\ref{cw4}) for $u$ has a mixed type since it also involves
   the chemical potential $\mu$. It turns out that for the
corresponding elliptic operator, it yields a non-homogeneous
boundary condition.  The \L ojasiewicz-Simon inequality for
homogeneous boundary
  conditions in the literature fails to apply. As a result, a non-trivial modification
  is required to treat the present problem. We succeed in deriving an extended \L ojasiewicz-Simon type
inequality involving boundary term, with which we able to show the
convergence result (for other applications, see \cite{WZ1,WZ2,Wu05,
WGZ2}). Besides, by delicate energy estimates and constructing
proper differential inequalities, we are able to obtain the
estimates for the convergence rate (in higher order norm). This in
some sense improves the previous result in the literature (see for
instance \cite[Theorem 1.1]{HJ01}) and can
 apply to other evolution equations (ref. \cite{Wu05,WGZ1}).

 \par The rest part of this paper is organized as follows: In Section 2 we introduce the functional settings and present
 some known results on existence and uniqueness of global solution and uniform compactness obtained in \cite{Gal,Gal2}.
 In Section3 we study the stationary problem. Section 4 is devoted to prove an extended \L ojasiewicz-Simon inequality with boundary term.
 In the final Section 5 we give the detailed proof of Theorem 1.1.


\section{Preliminaries}
\setcounter{equation}{0} We shall use the functional settings
introduced in
\cite{Gal,Gal2}.\\
For $u\in C(\overline{\Omega })$, we identify $u$ with the vector
$U=(u\big|_{\Omega },u\big|_{\Gamma })\in C(\Omega )\times C(\Gamma
)$. We define $\mathcal{H}=L^{2}(\Omega)\oplus L^{2}(\Gamma) $ to be
the completion of $C(\overline{\Omega })$ with respect to the
following norm
\begin{equation}
\|u\|_{\mathcal{H}}=\left(\|u\|^2+\|u\|^2_{L^2(\Gamma)}\right)^\frac12.
\label{cwe2.1}
\end{equation}
For any $g\in \mathcal{H}$, consider the elliptic boundary value
problem
 \be \left\{\begin{array}{l}  -\Delta u =g, \quad\text{in }\Omega ,\\
\partial _{\nu}u + u=g, \quad\text{on }\Gamma.\label{cwEWE}
  \end{array}
  \right.
 \ee
We can associate it with the following bilinear form on
$H^1(\Omega)$:
\begin{equation}
a(u,v)=\int_{\Omega }\nabla u \cdot \nabla u dx+\int_{\Gamma }u v
dS,
\label{cwe2.4}
\end{equation}
for all $u,v \in H^{1}(\Omega )$. Then it defines a strictly
positive self-adjoint unbounded operator $A: D(A)=\{u\in
H^1(\Omega)|Au\in \mathcal{H}\}\rightarrow \mathcal{H}$ such that
\begin{equation}
\langle Au, v \rangle _{\mathcal{H}} =a(u,v),\quad \forall\ u\in
D(A),\ v\in H^1(\Omega). \label{cwe2.5}
\end{equation}
Then by Lax-Milgram theorem, it follows that the operator $A$ is a
bijection from $D(A)$ into $\mathcal{H}$ and $A^{-1}:\mathcal{H} \to
\mathcal{H}$ is a linear, self-adjoint and compact operator on
$\mathcal{H}$ (see \cite[Section 2]{Gal2} or \cite[Section 4]{Gal}).
In other words, for any $ g\in \mathcal{H}$, $A^{-1}g$ is the unique
solution to \eqref{cwEWE}.

We can then consider the \textit{weak energy space} $\mathrm{X}$
endowed with the following norm,
\begin{equation} \|v\|_{\mathrm{X}}^{2}=\|
A^{-1/2}v\|_{\mathcal{H}}^{2}=\langle A^{-1}v,v\rangle
_{\mathcal{H}},\qquad \forall\  v\in \mathcal{H}.\label{cwe2.8}
\end{equation}
It follows that
 \be \big\langle
u,v\big\rangle _{\mathrm{X}}=\big\langle u,A^{-1}v \big\rangle
_{\mathcal{H}}, \quad \forall \ u\in H^1(\Omega),\  v \in
\mathrm{X}. \label{cwe2.10} \ee
 In particular, for all $v \in \mathrm{X}$ and $u=A^{-1}v$
  \be
\|v\|_{\mathrm{X}}^{2}=\langle A^{-1}v,v\rangle _{\mathcal{H}}=
\langle u, Au \rangle _{\mathcal{H}}=a(u,u). \label{cwe2.10a} \ee
For more detailed discussions, we refer to \cite{Gal,Gal2}.

The existence and uniqueness of global solution to
\eqref{cw1}-\eqref{cw5} has been obtained in \cite{Gal,Gal2}. The
results in \cite[Section 3,4]{Gal2} and \cite[Section 4]{Gal} in
particular imply that
 \bt Let
(\textbf{F1})--(\textbf{F3}) be satisfied. For any initial datum
$u_0\in V$, problem \eqref{cw1}-\eqref{cw5} admits a unique global
solution $u(t,x)$ which defines a
 global semiflow on $V$. Moreover, $u(t,x)$ belongs to $C^\infty$ for
 $t>0$.
\et
 \noindent The total free energy
 \be
 E(u)= \int_{\Omega} \left(\frac{1}{2} |\nabla u|^2 + F(u)\right) dx
 + \int_\Gamma\left(\frac{1}{2}\left|\nabla_\parallel u \right|^2 +\frac{1}{2} u^2\right) dS. \label{2.1}\ee
where $F(s)=\int_0^s f(z)dz$, serves as a Lyapunov functional for
problem (\ref{cw1})--(\ref{cw5}). In other words, for the smooth
solution $u$ to problem (\ref{cw1})--(\ref{cw5}),
 we have
  \be     \frac{d}{dt} E(u) +
\int_{\Omega} |\nabla \mu |^2 +\int_{\Gamma} |\mu|^2 dS = 0.
\label{cw2.2}\ee

Uniform bounds for the solution which yield the relative compactness
in $H^3(\Omega)\cap H^3(\Gamma)$ can be seen from \cite[Proposition
3.3, Theorem 3.5]{Gal2}, here we state the result without proof.
\begin{lemma} \label{cwcomp}
Let (\textbf{F1})--(\textbf{F3}) hold and $\gamma \in [0,1/2)$.
Then, for any initial datum $u_0\in V$, the solution of
\eqref{cw1}--\eqref{cw5} satisfies the following dissipative
estimates, namely, for any $\delta>0$, there hold
\begin{gather}
\|u(t)\|_{H^{3+\gamma}(\Omega)}^{2}+\|u(t)\|_{H^{3+\gamma}(\Gamma)}^{2}\leq
C_\delta, \quad t\geq \delta>0, \label{cwe3.23}
\end{gather}
and (ref. \cite[(3.34)]{Gal2}),
\begin{gather}
\|u_t\|_{H^1(\Omega)}^{2}+\|u_t\|_{H^1(\Gamma)}^{2}\leq C_\delta,
\quad t\geq \delta>0, \label{cwe3.23a}
\end{gather}
where $C_\delta>0$ depends only on $\|u_0\|_V$ and $\delta$.
\end{lemma}

 For any initial datum $u_0\in V$, the
$\omega$-limit set of $u_0$ is defined as follows:
$$\omega(u_0) = \{ \psi(x) \mid\ \exists \ \{t_n\}\  \text{such that}\
u(t_n,x) \rightarrow \psi(x)\in\ V,\ \text{as}\  t_n\rightarrow
+\infty \}.$$ Then we have
 \bl
 \label{cwlas} For any $u_0 \in V$, the
$\omega$-limit set of $u_0$ is a compact
connected subset in $H^3(\Omega)\cap H^3(\Gamma)$. Furthermore, \\
(i) $\omega(u_0)$ is  invariant under the nonlinear semigroup $S(t)$
defined by the solution $u(x,t)$, i.e, $S(t)\omega(u_0) =
\omega(u_0)$ for
all $ t \geq 0$. \\
(ii) $E(u)$ is constant on $\omega(u_0)$. Moreover, $\omega(u_0)$
consists of equilibria.
 \el
\begin{proof}
Since our system has a continuous Lyapunov functional $E(u)$, the
conclusion of the present lemma follows from Lemma \ref{cwcomp} and
the well-known results in the dynamical system (e.g. \cite[Lemma~
I.1.1]{T}).
 Thus, the lemma is proved.
\end{proof}

\section{Stationary Problem}
\setcounter{equation}{0} In this section we study the stationary
problem. The stationary problem corresponding to
\eqref{cw1}--\eqref{cw5} is \be \left\{\begin{array}{l}
  \Delta\widetilde{\mu}=0,\qquad  x \in \Omega,\\
 - \Delta \psi +f(\psi)=\widetilde{\mu}, \quad  x \in \Omega,\\
 \partial_\nu \widetilde{\mu}+\widetilde{\mu}=0, \qquad x\in \Gamma,\\
  -\Delta_\parallel \psi+\partial _\nu \psi
+\psi =\widetilde{\mu}, \quad x\in \Gamma.
  \end{array}
 \label{cwsta1}
 \right.
 \ee
Then it immediately follows that $\widetilde{\mu}=0$ and the
stationary problem is reduced to \eqref{cwsta}.

\bl Let (\textbf{F1})--(\textbf{F3}) be satisfied. Suppose that
$\psi \in H^3(\Omega)\cap H^3(\Gamma)$ satisfies (\ref{cwsta}). Then
$\psi$ is a critical point of the functional $E(u)$ over $V$.
Conversely, if $\psi \in V$ is a critical point of $E(u)$, then
$\psi \in C^{\infty}$ and it is a classical solution to problem
(\ref{cwsta}). \el
 \begin{proof}
 The proof is similar to \cite[Lemma 2.1]{WZ1}. The $C^{\infty}$ regularity for
 solution $\psi$ follows from the elliptic regularity
 for \eqref{cwsta} (see e.g., \cite[Corollary A.1, Lemma A.2]{MZ05})
 and a bootstrap argument.
 \end{proof}
 \bl Let (\textbf{F1})--(\textbf{F3}) be satisfied. The functional $E(u)$ has at least a
    minimizer $v\in V$ such that
    \be
    E(v)= \dis{\inf_{u\in V}}E(u).\ee
    In other words, problem (\ref{cwsta}) admits at least a
    classical solution. \el
\begin{proof}
 From assumption (\textbf{F3}), there exists $\delta_0>0$ such that
  \be  \liminf_{|s|\rightarrow +\infty} f'(s)\geq \delta_0.\ee
 Then there exists $N_1=N_1(\delta_0)>0$ such that
\be   f'(s)\geq \frac12 \delta_0,\qquad |s|\geq N_1.\ee It follows
that \be \liminf_{s\rightarrow +\infty} f(s)\geq 1,\quad
\liminf_{s\rightarrow -\infty} f(s)\leq -1.\label{cwff}\ee Since
$F'(s)=f(s)$, then we can deduce from \eqref{cwff} that \be
\liminf_{|s|\rightarrow +\infty} F(s)\geq 1.\ee Therefore, there
exists $N_2\geq 0$ such that
 \be F(s)\geq 0,\qquad |s|\geq N_2.\ee
 This indicates that
 \be \int_\Omega F(u)dx =\int_{|u|> N_2} F(u)dx+ \int_{|u|\leq N_2}
 F(u)dx\geq |\Omega|\min_{|s|\leq N_2} F(s)>-\infty.\ee
$E(u)$ can be written in the form:
    \be
    E(u)=\frac{1}{2}\| u\|_V^2+\mathcal{F}(u)\ee
    with
    \be
    \mathcal{F}(u)=\int_{\Omega}F(u)dx.\ee
It follows that $E(u)$ is bounded from below on $V$, namely,
 \be E(u)\geq \frac12\|u\|^2_V+C_f,\label{bl}\ee
where $C_f:=|\Omega|\min_{|s|\leq N_2} F(s)$. It's easy to see that
constant $C_f$ depends only on $f$ and $\Omega$. Therefore, there is
a minimizing sequence $u_n\in V$ such that
     \be E(u_n) \rightarrow \dis{\inf_{u\in V}}E(u).\label{cw3.3}\ee
       It follows from \eqref{bl}
       that $u_n$ is bounded in $V$. It turns out from the weak compactness
       that there is a subsequence, still denoted by $u_n$,
       such that $u_n$ weakly converges to $v$ in $V$. Thus,
       $v \in V$. We infer from the Sobolev imbedding theorem that the
       imbedding $V\subset H^1(\Omega)\hookrightarrow  L^\gamma(\Omega) $ ( $1\leq \gamma< \frac{n+2}{n-2}$)
        is compact. As a result, $u_n$ strongly converges to $v$ in $L^\gamma(\Omega)$. It turns out
        from the assumption (\textbf{F2}) that $\mathcal{F}(u_n) \rightarrow
    \mathcal{F}(v)$.
    Since $\|u\|^2_{V}$ is weakly lower semi-continuous, it follows
     from (\ref{cw3.3}) that $E(v)=\dis{\inf_{u\in V}}E(u)$.

The proof is completed.
\end{proof}

\section{Extended \L ojasiewicz-Simon Inequality }
 \setcounter{equation}{0}

In what follows, we prove a suitable version of extended ~\L
ojasiewicz-Simon inequality required in the proof of our main
result.

Let $\psi$ be a critical point of $E(u)$. We consider the following
linearized operator
 \be  L(v) h  \equiv -\Delta h+ f'(v+\psi)h\label{cw3.1} \ee
with the domain being defined as follows.
   \be Dom(L(v)) = \{h \in H^2(\Omega)\cap H^2(\Gamma) :  -
   \Delta_\parallel h + \partial_\nu h +h \mid_\Gamma =
   0\}:=\mathcal{D}.
   \ee
The equivalent norm on $\mathcal{D}$ is  \be
\|u\|_\mathcal{D}:=\|u\|_{H^2(\Omega)}+\|u\|_{H^2(\Gamma)}.\ee
 It's
obvious that $\mathcal{D}\subset L^2(\Omega)$ is dense in
$L^2(\Omega)$, and $L(v)$ maps $\mathcal{D}$ into $L^2(\Omega)$. In
analogy to \cite[Lemma 2.3]{WZ1}, we know that $L(v)$ is
self-adjoint.

Associated with $L(0)$, we define the bilinear form $b(w_1, w_2)$ on
$V$ as follows. \be b(w_1, w_2)= \int_\Omega (\nabla w_1 \cdot
\nabla w_2 + f'(\psi)w_1w_2 )dx +
  \int_\Gamma \left( \nabla_\parallel w_1 \cdot \nabla_\parallel w_2 + w_1 w_2\right) dS
\label{cw2.8}\ee Then, the same as for the usual second order
elliptic operator, $L(0)+\lambda I$ with $\lambda>0$ being
sufficiently large is invertible and its inverse is compact in
$L^2(\Omega)$. It turns out from the Fredholm alternative theorem
that $Ker(L(0))$ is finite-dimensional. It is well
 known that
 \be Ran(L(0))=(Ker(L(0))^*)^\perp. \ee
Thus, we infer from the fact that $L(0)$ is a self-adjoint operator
that \be Ran(L(0))=(Ker(L(0)))^\perp, \qquad Ran(L(0))\oplus
Ker(L(0))=L^2(\Omega). \ee Next we introduce two orthogonal
projections $\Pi_K$ and $\Pi_R$ in $L^2(\Omega)$, namely, $\Pi_K$ is
the projection onto the kernel of $L(0)$ while $\Pi_R$ is the
projection onto the range of $L(0)$.   Then we have the following
result. \bl For
$$L(0)w=f_R$$
with $f_R\in L^2(\Omega)$,  there exists a
 unique solution $w_R\in \mathcal{D}$ and the following estimate holds:\\
 \be
 \| w_R\|_{\mathcal{D}}\ \leq C
 \|f_R\|. \label{cw2.9}\ee
\el
\begin{proof}
By the Fredholm alternative theorem and the regularity theorem for
the elliptic operator (see \cite{MZ05}), we have a function $w \in
\mathcal{D}$ such that $L(0)w=f_R$. Moreover $w$ is unique if we
require $w\in (KerL(0))^\perp$, and (\ref{cw2.9}) follows from the
elliptic regularity theory.
\end{proof}
\noindent Let $\mathcal{L}(v): \mathcal{D}\rightarrow L^2(\Omega)$
be defined as follows:
 \be \mathcal{L}(v) w = \Pi_K w + L(v) w.
\label{cw3.17}\ee Then it follows from the above lemma that
$\mathcal{L}(0)$ is bijective and its inverse $\mathcal{L}^{-1}(0)$
is a bounded linear operator from $ L^2(\Omega)$ to $ \mathcal{D}$.
\bl\label{cwlemma1}
 There exists a small constant $ \beta <1$ depending on $\psi$ such that for any $v\in \mathcal{D}$,
$\| v\|_{H^2(\Omega)} \leq \beta$ and $f \in L^2(\Omega) $, \be
 \mathcal{L}(v) w = f \label{cw3.24}
 \ee
admits a unique solution $w$ such that $w\in  \mathcal{D}$ and the
following estimate holds, \be \| w \|_{\mathcal{D}} \leq C\| f \|.
\label{cw3.25aa} \ee \el
\begin{proof}
It follows from the above lemma that $\mathcal{L}(0)$ is bijective
and its inverse $\mathcal{L}^{-1}(0)$ is a bounded linear operator
from $ L^2(\Omega)$ to $ \mathcal{D}$. We rewrite (\ref{cw3.24})
into the following form: \be   (\mathcal{L}^{-1}(0)(
\mathcal{L}(v)-\mathcal{L}(0))+I) w =\mathcal{L}^{-1}(0) f.
\label{cw3.24a} \ee From the definition, we have
~$(\mathcal{L}(v)-\mathcal{L}(0))w=(f'(v+\psi)-f'(\psi))w$.

We infer from Sobolev imbedding theorem that for any
$\|v\|_{H^2}\leq \beta\ll 1$, there holds  \be \|
(f'(v+\psi)-f'(\psi))w\|\leq C\|v
\|_{H^2(\Omega)}\|w\|_{\mathcal{D}}.\ee Therefore, it follows that
when $\beta$ is sufficiently small,
 $\mathcal{L}^{-1}(0)( \mathcal{L}(v)-\mathcal{L}(0))$ is a
 contraction from $\mathcal{D}$ to $\mathcal{D}$:
  \be
  \|\mathcal{L}^{-1}(0)(
  \mathcal{L}(v)-\mathcal{L}(0))\|_{L(\mathcal{D}, \mathcal{D})}\leq
  \frac{1}{2}.\ee
By the contraction mapping theorem, (\ref{cw3.24a}) is uniquely
  solvable which implies that when $\|v\|_{H^2(\Omega)}\leq \beta$, $\mathcal{L}(v)$ is invertible, and (\ref{cw3.25aa}) holds.

  Thus, the lemma is proved.
 \end{proof}
\noindent Let $\psi$ be a critical point of $E(u)$.
Denote~$u=v+\psi$ and
 \be
\mathcal{E}(v) = E(u) = E(v+\psi).\label{cw3.aa}\ee
  Let
  \be M(v)=-\Delta(v + \psi) + f(v+\psi). \ee
  Then for any $v\in \mathcal{D}$, $M(v)\in
  L^2(\Omega)$.

First, we prove the following \L ojasiewicz-Simon inequality for the
homogeneous boundary condition corresponding to the nonhomogeneous
one \eqref{cw4}.
 \bl \label{cwls1} Let
$\psi$ be a critical point of $E(u)$. There exist constants
$\theta^*\in (0,\frac{1}{2})$ and $\beta^*\in(0,\beta)$ depending on
$\psi$ such that for any $ w \in \mathcal{D}$, if
$\|w\|_\mathcal{D}<\beta ^* $, there holds
   \be
   \| M(w) \| \ \geq\  | \mathcal{E}(w) - E(\psi)|^{1-\theta^*}.\label{cw3.16a}\ee
   \el
\begin{proof}
 Let $ \mathcal{N}: \mathcal{D} \mapsto L^2(\Omega)$ be the
 nonlinear operator defined as follows
 \be \mathcal{N}(w)=\Pi_K w+ M(w).\ee Then $\mathcal{N}(w)$ is differentiable and
 \be
D\mathcal{N}(w)h=\mathcal{L}(w)h. \ee

By the result in \cite{MR}, we know that
 \bl The mapping $
L^\infty(\Omega)\ni u \rightarrow f(u)\in L^{\infty}(\Omega)$ is
analytic. \el
 It easily follows from Lemma~\ref{cwcomp}, Sobolev
imbedding theorem and above lemma that $\mathcal{N}(w)$ is analytic.
Since $\mathcal{L}(0)$ is invertible, by the abstract implicit
function theorem (for the analytic version see e.g. \cite[Corallary
~4.37, p.172]{Zei}), there exist neighborhoods of the origin $W_1(0)
\subset \mathcal{D}$, $W_2(0)\subset L^2(\Omega)$ and an analytic
inverse mapping $\Psi$ of $\mathcal{N}$ such that
$\Psi:W_2(0)\rightarrow W_1(0)$
 is~1-1 and onto. Besides,
 \bea && \mathcal{N}(\Psi(g)) = g \qquad
\forall g \in W_2(0),
   \label{3.20}\\
  & & \Psi(\mathcal{N}(v)) = v \qquad  \forall v \in W_1(0),\label{cw3.21}
\eea and in analogy to the argument in \cite[Lemma 1, pp.75]{S96}
(see also \cite[Lemma 5.4]{J981}) we can show that
  \bea && \| \Psi(g_1) - \Psi(g_2) \|_{\mathcal{D}}\ \leq C
\|g_1 - g_2 \| \qquad  \forall g_1, g_2 \in W_2(0),\label{cw3.22}\\
  & & \|\mathcal{N}(v_1) - \mathcal{N}(v_2) \|\
   \leq C \|  v_1-v_2 \| _{\mathcal{D}} \qquad \forall v_1, v_2
   \in W_1(0). \label{cw3.22a}
\eea
 Let
$\phi_1,...,\phi_m$ be the orthogonal unit vectors spanning
$Ker(L(0))$.
\\Since $\Psi$ is analytic, it turns out that \be
\Gamma(\xi):=\mathcal{E}\left(\Psi\left(\sum^m_{i=1}\xi_i\phi_i\right)\right)
\label{cw3.30} \ee is analytic with respect to
$\xi=(\xi_1,...,\xi_m)$ with $|\xi|$ sufficiently small  such that
$\Pi_K w=\dis{\sum^{m}_{i=1}}\xi_i\phi_i\in W_2(0)$.

With the aid of $\Gamma(\xi)$ which is an analytic function defined
in $\mathbb{R}^m$, we are able to apply the \L ojasiewicz
inequality. By the standard argument (see e.g. \cite{J981}), we can
show that, there exist constants $\theta^*\in (0,\frac{1}{2})$ and
$\beta^*\in(0,\beta)$ depending on $\psi$ such that for any $w\in
\mathcal{D}$ with $\|w\|_{\mathcal{D}}<\beta^*$, there holds \be
   \|  M(w) \| \ \geq\  | \mathcal{E}(w) - \Gamma(0)|^{1-\theta^*},\ee
which is exactly (\ref{cw3.16a}). The details are omitted.
\end{proof}

 Now we are in a position to prove the following extended \L ojasiewicz-Simon
 inequality with boundary term.
\bl \label{cwls} Let $\psi$ be a critical point of $E(u)$. Then
there exist constants $\theta\in(0,\frac{1}{2})$, $\beta _0\in
(0,\beta)$ depending on $\psi$
   such that for any $ u \in H^3(\Omega)$, if $\|
   u-\psi\|_{H^2(\Omega)}< \beta _0$, the following inequality holds,
   \be
   \| M(v) \| + \| -\Delta_\parallel u + \partial_\nu u + u\|_{L^2(\Gamma)}
   \ \geq\  | E(u) - E(\psi)|^{1-\theta}.\label{cw3.16}\ee
   \el
\begin{proof} For any  $ u \in H^3(\Omega)$, let $v=u-\psi$. Then $v\in H^3(\Omega)$. \\
We consider the
following elliptic boundary value problem:
\be \left\{\begin{array}{l}  - \Delta w=-\Delta v,\qquad x \in \Omega,\\
-\Delta _\parallel w +\partial_\nu w + w=0, \qquad x\in \Gamma.
  \end{array}
 \label{cwaux}
 \right.
 \ee
Since $\Delta v\in L^2(\Omega)$, similar to the previous discussion
for $L(0)$, it follows that equation \eqref{cwaux} admits a unique
solution $w\in \mathcal{D}$. From the $H^2$-regularity for
\eqref{cwaux} (see e.g., \cite[Appendix Lemma A.1]{MZ05}), it turns
out that
 \be
\|w\|_{H^2(\Omega)}+\|w\|_{H^2(\Gamma)}\leq C\|\Delta v\|\leq
C\|v\|_{H^2(\Omega)}.\label{cwauxe3}\ee
 Hence, there
exists $\widetilde{\beta}\in (0,\beta)$ such that for
$\|v\|_{H^2(\Omega)}<\widetilde{\beta}$ we have  \be
\|w\|_\mathcal{D}< \beta^{*}.\ee Here $\beta^{*}$ is the constant in
Lemma~\ref{cwls1}. Thus, \eqref{cw3.16a} holds for $w$.

On the other hand, \eqref{cwaux} can be rewritten in the following
form
\be \left\{\begin{array}{l}  - \Delta (w-v)=0,\qquad x \in \Omega,\\
- \Delta _\parallel (w-v) +\partial_\nu (w-v) + (w-v)= \Delta
_\parallel v-\partial_\nu v - v,\quad x \in \Gamma.
  \end{array}
 \label{cwaux1}
 \right.
 \ee
Again from \cite[Appendix Lemma A.1]{MZ05}, there holds
 \be
\|w-v\|_{H^1(\Omega)}+\|w-v\|_{H^1(\Gamma)}\leq C\| \Delta
_\parallel v-\partial_\nu v - v\|_{L^2(\Gamma)}.\label{cwesh1}\ee

\noindent By straightforward computation,  \bea \| M(w)\|&\leq&
\left(\| M(v)\|+C
\|v-w\|_{H^1(\Omega)}\right)\non\\
&\leq& \left(\| M(v)\|+C \| \Delta _\parallel v-\partial_\nu v -
v\|_{L^2(\Gamma)}\right). \label{cwvw1}\eea Meanwhile, it follows
from Newton-Leibniz formula that
\bea&& \mid E(w+\psi)- E(v+\psi)\mid \non\\
 & \leq & \left|\int_0^1\int_\Omega  M( v + t(w - v)) (v-w)
dxdt\right|\non\\
&&\
 +\left|\int_0^1\int_\Gamma (1-t)(\Delta _\parallel v-\partial_\nu v
- v)(v-w)  dSdt\right|
 \non \\
& \leq & C \left(\| M(v)\|+\|\Delta _\parallel v-\partial_\nu v -
v\|_{L^2(\Gamma)}\right)\| \Delta
_\parallel v-\partial_\nu v - v\|_{L^2(\Gamma)}\non \\
& \leq & C \left(\| M(v)\|+\|\Delta _\parallel v-\partial_\nu v -
v\|_{L^2(\Gamma)}\right)^2. \label{cwvw2}\eea
Since \bea &&|E(w+\psi)- E(\psi)\mid^{1-\theta^*}\non\\
 &\geq&   \mid
E(v+\psi)- E(\psi)\mid^{1-\theta^*}- \mid E(w+\psi)-
E(v+\psi)\mid^{1-\theta^*} \label{cwvw3},\eea and $0 < \theta^*<
\frac{1}{2}$,  $ 2(1-\theta^*)-1>0$, then we infer from
(\ref{cwvw1})--(\ref{cwvw3}) that
 \be C(\| M(v) \|+\|  \Delta
_\parallel v-\partial_\nu v -  v\|_{L^2(\Gamma)})\ \geq
 |E(u) - E(\psi) |
^{1-\theta^*}. \non\ee Taking $\varepsilon\in (0,\theta^*)$ and
$\beta_0\in (0,\widetilde{\beta})$, such that for $\|v\|_{H^2}<
\beta_0$,
 \be \frac{1}{C} \mid E(v+ \psi) -
E(\psi)|^{-\varepsilon}\geq 1. \label{cw4.ddd}\ee Let
$\theta=\theta^* -\varepsilon\in(0,\frac{1}{2})$, then for
$\|v\|_{H^2} < \beta_0$, there holds \be \| M(v)\|+\|\Delta
_\parallel v-\partial_\nu v - v\|_{L^2(\Gamma)} \ \geq
 | E(u) - E(\psi) |
^{1-\theta}, \label{cw4.30}\ee which is exactly \eqref{cw3.16} by
the definition of $v$.
\end{proof}


\section{Convergence to equilibrium and convergence rate}
\setcounter{equation}{0}

After the previous preparations, we now proceed to finish the proof
of Theorem \ref{cwcon}.

\textbf{Part I. Convergence to Equilibrium} \\
From the previous results, there exists an increasing sequence
$\{t_n\}_{n\in\mathbb{N}}, \, t_n\rightarrow
   +\infty$ and $\psi\in \omega(u_0)$ such that
   \be \lim_{t_n\rightarrow +\infty} \|u(t_n,x)- \psi(x)\|_{H^3(\Omega)}=0. \label{cw3.ee}\ee
On the other hand, it follows from (\ref{cw2.2}) that $E(u)$ is
   decreasing in time. We now consider all possibilities.

   (1). If there is a $t_0>0$ such that at this time
   $E(u)=E(\psi)$, then for all $t>t_0$, we deduce from
   (\ref{cw2.2}) that $\|\mu(t)\|_{H^1(\Omega)}\equiv 0$.
On the other hand,    it follows from (\ref{cw1}) (\ref{cw3}) that
\be \left\{\begin{array}{l}
-\Delta \mu =-u_t, \quad  x \in \Omega,\\
 \partial _\nu \mu+ \mu=-u_t, \quad  x \in \Gamma.   \label{cw3.ff}   \end{array}  \right.
\ee Then by \eqref{cwe2.10a} we have
 \be
\|u_t\|^2_{\mathrm{X}}=\int_\Omega \mu u_tdx +\int_\Gamma \mu u_t
dS=\|\mu\|^2_{H^1(\Omega)}.\label{cweqnorm}\ee This implies that
$\|u_t\|_\mathrm{X}\equiv 0$, i.e., $u$ is independent of $t$ for
all $t>t_0$. Since $u(x, t_n)\rightarrow \psi$, then (\ref{cwconr})
holds.

   (2). If for all $t>0$, $E(u)>E(\psi)$, and there is $t_0>0$ such that for all $t\geq t_0$,
   $v=u-\psi$ satisfies the condition of Lemma \ref{cwls}, i.e., $\|u-\psi\|_{H^2(\Omega)} < \beta_0$,
   then for the constant $\theta\in(0,\frac12)$ in Lemma  \ref{cwls}, we have
    \be
   -\frac{d}{dt}(E(u)-E(\psi))^{\theta}=-\theta
   (E(u)-E(\psi))^{\theta-1}\frac{d E(u)}{dt}.\label{cw3.51}\ee
   From \eqref{cw2}, $M(v)=\mu$. Then it follows from \eqref{cwH1}(\ref{cw2.2}) and Lemma \ref{cwls} that
\be -\frac{d}{dt}(E(u) - E(\psi))^\theta\geq \theta \frac{\|\nabla
\mu\|^2+\|\mu\|^2_{L^2(\Gamma)}}{\|\mu\|+\|\mu\|_{L^2(\Gamma)}} \geq
C_\theta \| \mu \|_{H^1(\Omega)}. \label{cw3.52}\ee Integrating from
$t_0$ to~$t$, \be (E(u) - E(\psi))^{\theta} + C_\theta \int_{t_0}^t
\| \mu \|_{H^1(\Omega)} d\tau\leq (E(u(t_0)) - E(\psi))^\theta.
\label{cw3.53}\ee Since, $E(u(t)) - E(\psi) \geq 0$, we have \be
\int_{t_0}^t \| \mu \|_{H^1(\Omega)} d\tau < +\infty,\qquad \forall\
t\geq t_0. \label{cw3.54}\ee
 Thus, \eqref{cweqnorm}(\ref{cw3.54}) imply that for all $t\geq t_0$,  \be \int^t
_{t_0}\|u_t\|_{\mathrm{X}}d \tau < +\infty, \label{cw3.55}\ee which
easily yields that as $t\rightarrow +\infty$, $u(t,x)$ converges in
$\mathrm{X}$. Since the orbit is compact in $H^3(\Omega)\cap
H^3(\Gamma)$, we can deduce from uniqueness of limit that
(\ref{cwconr}) holds.

 (3).  It follows from (\ref{cw3.ee}) that for  any $ \varepsilon\in (0,\beta_0)$,
there exists  $N \in \mathbb{N}$ such that when $n \geq N$, \bea &&
\parallel
 u(t_n,\cdot) - \psi
 \parallel_{\mathrm{X}}\  \leq\    \|u(t_n,\cdot)-\psi\|_{H^3(\Omega)}< \frac{\varepsilon}{2}, \label{3.55a}\\
 &&\frac{1}{C_\theta} (E(u(t_n)) - E(\psi))^\theta <
 \frac{\varepsilon}{2}. \label{cw3.56}\eea
Define
 \be \bar{t}_n = \sup \{\ t > t_n \mid \ \parallel
u(s,\cdot) - \psi \parallel_{H^2(\Omega)} < \beta _0,\ \forall s\in
[t_n , t] \}. \label{cw3.57}\ee (\ref{cw3.ee}) and continuity of the
orbit in $H^2(\Omega)$ yield that $\bar{t}_n >t_n$ for all $n \geq
N$.

Then there are two possibilities:\\
 (i). If there exists $n_0\geq N$
such that $\bar{t}_{n_0}=+\infty$, then from the previous
discussions in (1) and (2), \eqref{cwconr} holds.\\
 (ii) Otherwise, for
all $n\geq N$, we have $t_n < \bar{t}_n < +\infty$, and for all
$t\in [t_n, \bar{t}_n]$, $E(\psi)<E(u(t))$. Then from (\ref{cw3.53})
with $t_0$ being replaced by $t_n$, and $t$ being replaced by
$\bar{t}_n$
 we deduce that \be \int_{t_n}^{\bar{t}_n}
\parallel\! u_t \!\parallel_{\mathrm{X}} d\tau \leq C_\theta (E(u(t_n)) - E(\psi))^\theta <
\frac{\varepsilon}{2}. \label{cw3.59}\ee Thus we have \be
\parallel\! u(\bar{t}_n) - \psi \!
\parallel_{\mathrm{X}}\ \leq\ \ \parallel\! u(t_n) - \psi
\!\parallel_{\mathrm{X}} + \int_{t_n}^{\bar{t}_n}
\parallel\! u_t \!\parallel_{\mathrm{X}} d\tau < \varepsilon,
\label{cw3.60}\ee which implies that when $n\rightarrow +\infty $,
$$u(\bar{t}_n) \rightarrow \psi \qquad \text{in} \ \ \mathrm{X}.$$
Since $\bigcup_{t\geq\delta}u(t)$ is relatively compact in
$H^2(\Omega)$, there exists a subsequence of $\{u(\bar{t}_n)\}$,
still denoted  by $\{u(\bar{t}_n)\}$ converging to $\psi$ in
$H^2(\Omega)$.  Namely, when $n$ is sufficiently large, we have
$$\parallel\! u(\bar{t}_n) - \psi \!\parallel_{H^2(\Omega)}\ <
\beta_0,$$ which contradicts the definition of $\bar{t}_n$ that
$\|u(\bar{t}_n,\cdot)-\psi\|_{H^2(\Omega)}=\beta_0$.\\

 \textbf{Part II.
Convergence Rate}\\
For $t\geq t_0$, it follows from Lemma \ref{cwls} and (\ref{cw3.52})
that
 \be \frac{d}{dt}(E(u) - E(\psi))+ C\left(E(u) -
E(\psi)\right)^{2(1-\theta)}\ \leq\ 0. \label{cw3.52aa}\ee As a
result,
 \be E(u(t))-E(\psi) \leq  C (1+t)^{-1/(1-2\theta)},
\qquad\forall\,t\geq t_0. \ee Integrate \eqref{cw3.52} on
$(t,\infty)$, where $t\geq t_0$, then we have
 \be \int_{t}^{\infty}
\|u_t\|_{\mathrm{X}} d\tau
           \leq C (1+t)^{-\theta/(1-2\theta)}.
\ee By adjusting the constant $C$ properly, we obtain
 \be
    \|u(t)-\psi\|_{\mathrm{X}}\leq C
    (1+t)^{-\theta/(1-2\theta)}, \quad t\geq 0.\label{cwonlow}
\ee Based on this convergence rate we are able to get the same
estimate for convergence rate in higher order norm by energy
estimates and proper differential inequalities.\\
Next we proceed to estimate $\|u-\psi\|_V$.\\
It follows from \eqref{cw1}-\eqref{cw4} and the stationary problem
\eqref{cwsta} that
 \be \left\{
\begin{array}{l}
\displaystyle{\frac{d }{d t}}(u-\psi)=\Delta \mu, \\
 \mu = -\Delta(u-\psi)+f(u)-f(\psi),
 \end{array} \label{cw1.1m}
  \right.
 \ee
with the boundary condition
 \be \left\{\begin{array}{l}
 -\Delta _\parallel (u-\psi)+
\partial_\nu (u-\psi) + (u-\psi)= \mu\\
 (u-\psi)_t+\partial_\nu \mu+\mu=0. \end{array}
 \label{cw1.2m}
 \right.
 \ee
Using \eqref{cw1.1m}\eqref{cw1.2m}, we take the inner product in
$\mathcal{H}$ of $A^{-1}(u-\psi)_t$ with $(u-\psi)$ to obtain
 \bea
&& \frac12\frac{d}{dt}\|u-\psi\|^2_{\mathrm{X}}+\|\nabla u-\nabla
\psi\|^2+\int_\Omega(f(u)-f(\psi))(u-\psi)dx\non\\
&&\ \ +\left\|\nabla_\parallel
(u-\psi)\right\|_{L^2(\Gamma)}^2+\|u-\psi\|_{L^2(\Gamma)}^2\non\\
&=& 0.\label{cwonm1}
 \eea
On the other hand, by \eqref{cw1.1m}\eqref{cw1.2m} and taking the
inner product in $\mathcal{H}$ of $(u-\psi)_t$ with $\mu$, we have
 \bea &&
\frac{d}{dt}\left(\frac12\|\nabla u-\nabla \psi\|^2+ \int_\Omega
 F(u)dx-\int_\Omega f(\psi)udx+\frac{1}{2} \left\|\nabla_\parallel
(u-\psi)\right\|_{L^2(\Gamma)}^2\right.\non\\
&&\ \ \left. +\frac{1}{2} \|u-\psi\|_{L^2(\Gamma)}^2\right)+\|\nabla
\mu\|^2+\|\mu\|_{L^2(\Gamma)}^2\non\\
&=&0.\label{cwonm2}
 \eea
Adding \eqref{cwonm1}\eqref{cwonm2} together, we have
 \bea
 && \frac{d}{dt}\left(\frac12\|u-\psi\|^2_{\mathrm{X}}
 +\frac{1}{2}\|u-\psi\|^2_{L^2(\Gamma)}+\frac{1}{2} \left\|\nabla_\parallel
(u-\psi)\right\|_{L^2(\Gamma)}^2+\frac12\|\nabla u-\nabla \psi\|^2 \right.\non\\
&&\quad \quad \left.+ \int_\Omega
 F(u)dx-\int_\Omega F(\psi)dx+ \int_\Omega f(\psi)\psi dx- \int_\Omega f(\psi)u
 dx\right)\non\\
 &&\ \ +\|\nabla
u-\nabla \psi\|^2+\left\|\nabla_\parallel
(u-\psi)\right\|_{L^2(\Gamma)}^2+\|u-\psi\|_{L^2(\Gamma)}^2+\|\nabla
\mu\|^2+\|\mu\|_{L^2(\Gamma)}^2\non\\
&=& -\int_\Omega(f(u)-f(\psi))(u-\psi)dx.\label{cwonm3}
 \eea
In what follows, we shall use the uniform bounds obtained in Lemma
\ref{cwcomp}. Without loss of generality, we set $\delta=1$ in Lemma
\ref{cwcomp}.

The Newton-Leibniz formula
 \be F(u)=F(\psi) +
f(\psi)(u-\psi) +\int_0^1\int_0^1 f'(sz u+(1-sz)\psi)(u-\psi)^2
dsdz, \ee
 yields that
 \bea & & \left\vert\int_\Omega F(u)dx-\int_\Omega
        F(\psi)dx+\int_\Omega f(\psi)\psi
         dx-\int_\Omega f(\psi)u dx\right\vert\non\\
&  =  & \left\vert \int_\Omega \int_0^1\int_0^1
f'(sz u+(1-sz)\psi)(u-\psi)^2 ds dz dx \right\vert\non\\
& \leq &
          \max_{s, z\in [0,1]} \|f'(sz u+(1-sz)\psi)\|_{L^3}\|u-\psi\|^2_{L^3}\non\\
&\leq & C(\|\nabla u-\nabla\psi\|\|u-\psi\|+ \|u-\psi\|^2)\non\\
&\leq & \frac14\|\nabla u-\nabla\psi\|^2+
        C\|u-\psi\|^2,\quad t\geq 1. \label{cwascr10}\eea
  and in a similar way, we have
  \bea & & \left|\int_\Omega(f(u)-f(\psi))(u-\psi)dx\right|\non\\
&= & \left\vert \int_\Omega \int_0^1
f'(s u+(1-s)\psi)(u-\psi)^2 ds dx \right\vert\non\\
&\leq & \frac14\|\nabla u-\nabla\psi\|^2+
        C\|u-\psi\|^2,\quad t\geq 1. \label{cwascr10a}\eea
Let
 \bea y_1(t)&=& \frac12\|u-\psi\|^2_{\mathrm{X}}
 +\frac{1}{2}\|u-\psi\|^2_{L^2(\Gamma)}+\frac{1}{2} \left\|\nabla_\parallel
(u-\psi)\right\|_{L^2(\Gamma)}^2+\frac12\|\nabla u-\nabla \psi\|^2 \non\\
&&\quad \quad + \int_\Omega
 F(u)dx-\int_\Omega F(\psi)dx+ \int_\Omega f(\psi)\psi dx- \int_\Omega f(\psi)u
 dx
 \eea
\eqref{cwascr10} indicates that there exist constants $C_1,C_2>0$
such that
 \be y_1(t)\geq
C_1\|u-\psi\|^2_V-C_2\|u-\psi\|^2,\quad t \geq 1.\label{cwascr11}\ee
On the other hand,
 \bea
\|u-\psi\|^2 &\leq& C\|u-\psi\|_V\|u-\psi\|_{\mathrm{X}}\non\\
&\leq& \varepsilon
\|u-\psi\|^2_V+C(\varepsilon)\|u-\psi\|_{\mathrm{X}}^2.\label{cwascr12}
 \eea
From
\eqref{cwonlow}\eqref{cwonm3}\eqref{cwascr10}\eqref{cwascr11}\eqref{cwascr12},
after taking $\varepsilon>0$ sufficiently small, we can see that
there exists a constant $\gamma>0$ such that
 \be
  \frac{d}{dt}y_1(t)+\gamma y_1(t)\leq
  C\|u-\psi\|_{\mathrm{X}}^2\leq C(1+t)^{-2\theta/(1-2\theta)}, \qquad t\geq 1.\label{cwonm4}
 \ee
As a result, \bea y_1(t)
 &\leq& y_1(1)
e^{\gamma (1-t)}+Ce^{-\gamma t}\int_1^t
(1+\tau )^{-2\theta/(1-2\theta)}d\tau\non\\
&\leq& Ce^{-\gamma t}+ Ce^{-\gamma
t}\int_0^t (1+\tau )^{-2\theta/(1-2\theta)}d\tau\non\\
&\leq& Ce^{-\gamma t}+Ce^{-\gamma t}
\left(\int_0^{\frac{t}{2}}e^{\gamma\tau}(1+\tau)^{-2\theta/(1-2\theta)}
d
\tau+\int_{\frac{t}{2}}^te^{\gamma\tau}(1+\tau)^{-2\theta/(1-2\theta)}d
\tau\right)
\non\\
&\leq& C e^{-\gamma t}+Ce^{-\gamma t}\left(e^{\frac{\gamma
}{2}t}\int_0^{\frac{t}{2}}(1+\tau)^{-2\theta/(1-2\theta)}
d \tau+C(1+t)^{-2\theta/(1-2\theta)}e^{\gamma t}\right)\non\\
&\leq& C(1+t)^{-2\theta/(1-2\theta)},\quad t\geq
1.\label{cwscr12}\eea
 \eqref{cwascr11}\eqref{cwascr12}\eqref{cwscr12}
imply that
 \bea
C_1\|u-\psi\|^2_V&\leq& y_1(t)+C_2\|u-\psi\|^2\non\\
&\leq& y_1(t)+C_2\varepsilon \|u-\psi\|^2_V+
C_2C(\varepsilon)\|u-\psi\|_{\mathrm{X}}^2.\eea Taking
$\varepsilon>0$ sufficiently small, it follows from
\eqref{cwonlow}\eqref{cwscr12} that
 \be \|u-\psi\|_V\leq
 C(1+t)^{-\theta/(1-2\theta)},\qquad t\geq 1.\label{cwscr13}\ee
 By the $C^\infty$ regularity of the solution, we are able to get
 the estimate for convergence rate in higher order norm.\\
Differentiating \eqref{cw1}--\eqref{cw4} respect to time $t$
respectively, we have
 \be u_{tt}= \Delta\mu_t, \quad
x\in \Omega, \label{cwee2rac}\ee
 \be \mu_{t} = - \Delta u_{t} + f'(u)u_t,\quad
x\in \Omega, \label{cwee2ra}\ee
 \be u_{tt}+\partial_\nu \mu_t+\mu_t=0, \quad x\in \Gamma,   \label{cwee2rad}
 \ee
 \be \mu_t=-\Delta _\parallel u_t+\partial_\nu u_t
+ u_t,\quad x\in \Gamma.   \label{cwee2rab}\ee
 Multiplying
\eqref{cwee2ra} by $u_t$ and integrating by parts on $\Omega$, using
\eqref{cw3}\eqref{cwee2rab}, we get
 \bea
&&\frac{1}{2}\frac{d}{dt}\left(\|\nabla
\mu\|^2+\|\mu\|^2_{L^2(\Gamma)}\right)+\|\nabla u_t\|^2+\int_\Gamma
\left(\left|\nabla_\parallel u_t\right|^2+u_t^2\right)dS\non\\
&=& -\int_\Omega f'(u)u_t^2dx. \eea
 Assumption
(\textbf{F3}) yields that there is a certain positive constant
$M_f\geq1$ such that \be f'(s)\geq -M_f,\qquad s\in \mathbb{R}.\ee
 Thus, \be -\int_\Omega f'(u)u_t^2dx\leq M_f\|u_t\|^2.\label{cwesfm}\ee
 It follows
from \eqref{cw1} that
 \bea\|u_t\|^2&=&-\int_\Omega \nabla u_t\cdot\nabla
\mu dx-\int_\Gamma\mu u_t dS- \|u_t\|^2_{L^2(\Gamma)} \non\\
&\leq& \|\nabla u_t\|\|\nabla
\mu\|+\|\mu\|_{L^2(\Gamma)}\|u_t\|_{L^2(\Gamma)}\non\\
&\leq& \varepsilon\|\nabla
u_t\|^2+\varepsilon\|u_t\|^2_{L^2(\Gamma)}
+\frac{1}{4\varepsilon}\|\nabla
\mu\|^2+\frac{1}{4\varepsilon}\|\mu\|^2_{L^2(\Gamma)},\label{cwee7ab}\eea
In \eqref{cwee7ab}, taking \be \varepsilon=
\frac{1}{2M_f},\label{cwesfm1}
 \ee  it follows that
 \bea
&&\frac{1}{2}\frac{d}{dt}\left(\|\nabla
\mu\|^2+\|\mu\|^2_{L^2(\Gamma)}\right)+\frac12\|\nabla
u_t\|^2+\|\nabla_\parallel u_t\|_{L^2(\Gamma)}^2+\frac12 \|u_t\|^2_{L^2(\Gamma)}\non\\
&\leq& C\left(\|\nabla
\mu\|^2+\|\mu\|^2_{L^2(\Gamma)}\right).\label{cwscr20} \eea
Multiplying \eqref{cwee2ra} by $u_{tt}$ and integrating by parts on
$\Omega$, using \eqref{cwee2rac}\eqref{cwee2rad}\eqref{cwee2rab}, we
get
 \bea
&&\frac{1}{2}\frac{d}{dt}\left(\|\nabla u_t\|^2+\int_\Omega
f'(u)u_t^2dx+\|\nabla_\parallel u_t\|_{L^2(\Gamma)}^2+\|u_t\|^2_{L^2(\Gamma)}\right)\non\\
&&\ \ \ + \ \|\mu_t\|^2_{L^2(\Gamma)}+\|\nabla
\mu_t\|^2\non\\
&=& \frac12\int_\Omega f''(u)u_t^3dx.\label{cwee17}\eea
 By \eqref{cwee7ab} and Lemma \ref{cwcomp}, the
righthand side of \eqref{cwee17} can be estimated as follows
 \bea &&\left|\int_\Omega f''(u)u_t^3dx\right|\non\\
 &\leq& C(
|u|_{L^\infty})\|u_t\|^3_{L^3}\leq C\|u\|_{H^2}\left(\|\nabla
u_t\|^{\frac{3}{2}}\|u_t\|^{\frac{3}{2}}+\|u_t\|^3\right)\non\\
&\leq& \frac18\|\nabla
u_t\|^{2}+C\|u_t\|^6+C\|u_t\|^3\non\\
 &\leq&
\frac18\|\nabla u_t\|^{2}+C\|u_t\|^2\non\\
&\leq& \frac14\left(\|\nabla u_t\|^2+\|u_t\|^2_{L^2(\Gamma)}\right)
+C\|\nabla \mu\|^2+C\|\mu\|^2_{L^2(\Gamma)},\qquad t\geq
1.\label{cwee21}\eea
 Then
\eqref{cwee17} becomes
 \bea &&\frac{1}{2}\frac{d}{dt}\left(\|\nabla
u_t\|^2+\int_\Omega
f'(u)u_t^2dx+\|\nabla_\parallel u_t\|_{L^2(\Gamma)}^2+\|u_t\|^2_{L^2(\Gamma)}\right)\non\\
&&\ \ \ + \ \|\mu_t\|^2_{L^2(\Gamma)}+\|\nabla
\mu_t\|^2\non\\
&\leq& \frac18\left(\|\nabla u_t\|^2+\|u_t\|^2_{L^2(\Gamma)}\right)
+C\|\nabla \mu\|^2+C\|\mu\|^2_{L^2(\Gamma)},\qquad t\geq
1.\label{cwee22}\eea Multiplying \eqref{cwee22} by
$\varepsilon_1\in(0,1]$ and adding the resultant to \eqref{cwscr20},
we obtain
 \bea
&&\frac{1}{2}\frac{d}{dt}\left(\|\nabla
\mu\|^2+\|\mu\|^2_{L^2(\Gamma)}+\varepsilon_1\|\nabla
u_t\|^2+\varepsilon_1\int_\Omega
f'(u)u_t^2dx+\varepsilon_1\left\|\nabla_\parallel
u_t\right\|_{L^2(\Gamma)}^2\right.\non\\
&&\ \ \left.+ \varepsilon_1\|u_t\|^2_{L^2(\Gamma)}
\right)+\frac14\left(\|\nabla
u_t\|^2+\|u_t\|^2_{L^2(\Gamma)}\right)+\left\|\nabla_\parallel
u_t\right\|_{L^2(\Gamma)}^2\non\\
&&\ \ +\varepsilon_1\|\mu_t\|^2_{L^2(\Gamma)}+\varepsilon_1\|\nabla
\mu_t\|^2\non\\
&\leq& C^*\left(\|\nabla
\mu\|^2+\|\mu\|^2_{L^2(\Gamma)}\right),\qquad t\geq 1.\label{cwee23}
\eea
 Let
 \be y_2(t)=\|\nabla
\mu\|^2+\|\mu\|^2_{L^2(\Gamma)}+\varepsilon_1\|\nabla
u_t\|^2+\varepsilon_1\int_\Omega
f'(u)u_t^2dx+\varepsilon_1\left\|\nabla_\parallel
u_t\right\|_{L^2(\Gamma)}^2+\varepsilon_1
\|u_t\|^2_{L^2(\Gamma)}\label{cwee24}\ee It follows from Lemma
\ref{cwcomp} that
 \be y_2(t)\leq C,\qquad t\geq 1.\ee
Taking \be \varepsilon_1=\frac{1}{M_f^2},\ee we can deduce from
\eqref{cwesfm}--\eqref{cwesfm1} that
 \be y_2(t)\geq \frac12\left(\|\nabla
\mu\|^2+\|\mu\|^2_{L^2(\Gamma)}\right)+\varepsilon_1\left(
\frac12\|\nabla u_t\|^2+\left\|\nabla_\parallel
u_t\right\|_{L^2(\Gamma)}^2+\frac12
\|u_t\|^2_{L^2(\Gamma)}\right).\label{cwcy2g}\ee Now we take
$\kappa>0$ such that
 \be \kappa (1+C^*)\leq \frac12.\ee
Next, we multiply \eqref{cwee23} by $\kappa$ and add the resultant
to \eqref{cwonm3}, then
\eqref{cwonm4}\eqref{cwascr10a}\eqref{cwscr13} yield that there
exists a constant $\widetilde{\gamma}>0$ such that
 \be
  \frac{d}{dt}[y_1(t)+\kappa y_2(t)]+\widetilde{\gamma} [y_1(t)+\kappa
  y_2(t)]
  \leq C\|u-\psi\|^2\leq C(1+t)^{-2\theta/(1-2\theta)}, \quad t\geq 1.\label{cwonmy2}
 \ee
Similar to \eqref{cwscr12}, we have
 \be y_1(t)+\kappa y_2(t)\leq C(1+t)^{-2\theta/(1-2\theta)},\quad t\geq
1.\label{cwscry2}\ee Hence, from
\eqref{cwascr11}\eqref{cwscr12}\eqref{cwscr13}\eqref{cwscry2} we
know \be y_2(t)\leq C(1+t)^{-2\theta/(1-2\theta)},\quad t\geq
1,\label{cwscry2a}\ee which together with \eqref{cwcy2g} gives the
following
 \be \| \mu\|_{H^1(\Omega)}+\|u_t\|_{H^1(\Omega)} +\|u_t\|_{H^1(\Gamma)}\leq C(1+t)^{-\theta/(1-2\theta)},\quad t\geq
1.\label{cwcy2ga}\ee
  By the
elliptic estimate (see \cite[Corollary A.1]{MZ05}),
 \be
\|u-\psi\|_{H^3(\Omega)}+\|u-\psi\|_{H^3(\Gamma)}\leq C
\left(\|\mu\|_{H^1(\Omega)}+\|f(u)-f(\psi)\|_{H^1(\Omega)}+\|\mu\|_{H^1(\Gamma)}\right).\label{cwee35}\ee
Lemma \ref{cwcomp} and Sobolev imbedding theorem imply that
 \be \|f(u)-f(\psi)\|_{H^1(\Omega)}\leq C\|u-\psi\|_V,\qquad t\geq
1.\label{cee36}\ee On the other hand, from \eqref{cw3.ff}, the
elliptic regularity theory and Sobolev imbedding theorem, we have
 \be \|\mu\|_{H^1(\Gamma)}\leq C\|\mu\|_{H^2(\Omega)}\leq
 C\left(\|u_t\|+\|u_t\|_{H^\frac12(\Gamma)}\right)\leq
 C\|u_t\|_{H^1(\Omega)}.\label{cwee36}\ee
 As a result, we can conclude from \eqref{cwscr13},
\eqref{cwcy2ga}--\eqref{cwee36} that
 \be \|u-\psi\|_{H^3(\Omega)}+\|u-\psi\|_{H^3(\Gamma)} \leq C(1+t)^{-\theta/(1-2\theta)},\quad t\geq 1.\label{cwee30}\ee
Summing up, the proof of theorem \ref{cwcon} is completed.

\br Following the same method, we can continue to get estimates of
convergence rate in higher order norm.\er
 \br We notice that, in order to get the
convergence rate estimates \eqref{cwscr13} \eqref{cwcy2ga}
\eqref{cwee30}, we have to use the uniform bound for the solution in
higher order norm, e.g. Lemma \ref{cwcomp}, which is not valid for
$t=0$. Thus, the constant $C$ in
\eqref{cwscr13}\eqref{cwcy2ga}\eqref{cwee30} depends on $\delta$ in
Lemma \ref{cwcomp}. More precisely, for any $\delta>0$ we have
 \be \|u-\psi\|_{H^3(\Omega)}+\|u-\psi\|_{H^3(\Gamma)}+\|u_t\|_{V} \leq
C_\delta (1+t)^{-\theta/(1-2\theta)},\quad \forall\ t\geq \delta.\ee
the constant $C_\delta$ depends on $\|u_0\|_V$ and $\delta$.
Moreover, \be \lim_{\delta\rightarrow 0^+} C_\delta=+\infty.\ee \er

\noindent\textbf{Acknowledgements:} The author is indebted to Prof.
Songmu Zheng for his enthusiastic help and encouragement. This work
is partially supported by the NSF of China under grant No.~10631020
and Fudan Postgraduate Innovation Project under grant No.~CQH
5928003.


\end{document}